\documentclass{amsart}

\usepackage{amsmath,amsthm,amssymb,amsfonts}

\newcommand{\sgn}{\operatorname{sgn}}
\newcommand{\pfaff}{\mathop{\mathrm{pfaff}}}

\numberwithin{equation}{section} 
\theoremstyle{plain}
\newtheorem{theo+}           {Theorem}      [section]
\newtheorem{prop+}  [theo+]  {Proposition}
\newtheorem{coro+}  [theo+]  {Corollary}
\newtheorem{lemm+}  [theo+]  {Lemma}
\newtheorem{defi+}  [theo+]  {Definition}

\theoremstyle{definition}
\newtheorem{exam+}  [theo+]  {Example}
\newtheorem{rema+}  [theo+]  {Remark}

\newenvironment{theorem}{\begin{theo+}}{\end{theo+}}

\newenvironment{corollary}{\begin{coro+}}{\end{coro+}}
\newenvironment{lemma}{\begin{lemm+}}{\end{lemm+}}

\newenvironment{remark}{\begin{rema+}}{\end{rema+}}

\begin{document}

\baselineskip 18pt
\larger[2]
\title
[Sums of triangular numbers]
{Sums of triangular numbers\\ 
 from the Frobenius determinant} 
\author{Hjalmar Rosengren}
\address
{Department of Mathematical Sciences
\\ Chalmers University of Technology and G\"oteborg
 University\\SE-412~96 G\"oteborg, Sweden}
\email{hjalmar@math.chalmers.se}
\urladdr{http://www.math.chalmers.se/{\textasciitilde}hjalmar}

\thanks{Research  supported by the Swedish Science Research
Council (Vetenskapsr\aa det)}

\begin{abstract}

We show that the denominator formula for the strange series of affine
super\-algebras, conjectured by
Kac and Wakimoto and proved by Zagier, follows from a classical
determinant evaluation of Frobenius. As a limit case, we
obtain  exact formulas for the number of representations of an
arbitrary number as a sum of $4m^2/d$ triangles, whenever $d\mid 2m$,
and $4m(m+1)/d$ triangles, when $d\mid 2m$ or $d\mid 2m+2$. 
This extends  recent results of
Getz and Mahlburg, Milne, and Zagier. 
\end{abstract}

\maketitle        


\section{Introduction}  

To count the number of representations of an integer $n$ as a sum
of $k$ triangular numbers is a classical problem \cite{d}. We
will denote this number by $\triangle_k(n)$. The most 
 fundamental results are  
\begin{subequations}\label{l}
\begin{align}\label{t2}\triangle_2(n)&= \sum_{d\mid
    4n+1}(-1)^{\frac 12(d-1)},\\
\label{t4}  \triangle_4(n)&=\sum_{d\mid 2n+1}d,\\
\label{t8}  \triangle_8(n)&=\sum_{ d\mid n+1,\,
  (n+1)/d \text{ odd}}d^3
\end{align}
\end{subequations}
(where $d$ is assumed to be positive), which may be compared
 with  the sums of squares formulas
\begin{subequations}
\begin{align}\label{s2}\square_2(n)&=4\sum_{d\mid n,\, d \text{ odd}}(-1)^{\frac12(d-1)}
,\\
\label{s4} \square_4(n)&=8\sum_{d\mid n,\,4\nmid d}d,\\
\label{s8} \square_8(n)&=
16\sum_{d\mid n}(-1)^{n+d} d^3;\end{align}
\end{subequations}
see \eqref{tkd}  below for the exact conventions used for
defining $\triangle_k$. 
The identities \eqref{t4} and \eqref{t8} were found by Legendre,
 while \eqref{s4} and \eqref{s8} are due to Jacobi. As for 
 \eqref{s2}, it was first published by
 Legendre in 1798. 
It implies \eqref{t2}, in view of
$$4\left(\frac{k(k+1)}2+\frac{l(l+1)}2\right)+1=(k+l+1)^2+(k-l)^2.$$

There have been many attempts to find exact formulas for 
$\triangle_k$ and $\square_k$ when $k\neq 2, 4, 8$. Some ten years ago, a breakthrough was made by Kac and Wakimoto \cite{kw}, who showed
that many results of this type can be obtained from denominator
formulas for affine superalgebras. In particular, they conjectured a
denominator formula for the  strange series $Q(k)$
of affine superalgebras (see \cite{ka}),
and showed that it would imply new formulas for $\triangle_{4m^2}$ 
(when $k=2m-1$) and
$\triangle_{4m(m+1)}$ (when $k=2m$).  When
$m=1$ one  recovers Legendre's results for $4$ and $8$ triangles.

Roughly speaking, the Kac--Wakimoto formulas correspond to
replacing the divisor 
sums in \eqref{l} by sums over solutions to an equation of the form
$k_1l_1+\dots+k_ml_m=y$. We  recall the exact statements below,  see
\eqref{kmt} and \eqref{kmt2}. It may be instructive to give here the
case $m=2$ explicitly, namely,
\begin{align}\label{t16}\triangle_{16}(n)&=\frac 1{2^7\cdot
  3}\sum_{\substack{k_1l_1+k_2l_2=2n+4\\ k_i \text{ and } l_i \text{
  odd positive}}}k_1k_2(k_1^2-k_2^2)^2,\\
\label{t24}\triangle_{24}(n)&=\frac 1{2^4\cdot
  3^2}\sum_{\substack{k_1l_1+k_2l_2=n+3\\ k_i \text{ positive},\ l_i \text{
  odd positive}}}k_1^3k_2^3(k_1^2-k_2^2)^2.\end{align}
By symmetry, one may impose the condition $k_1<k_2$ if the
  right-hand sides are multiplied by $2$; this is of course preferable for the
  purpose of computation.

The denominator formula for the strange series was proved by Zagier
\cite{z}, using elliptic function identities. He also gave
a second proof of the corresponding triangular number identities,
using modular forms. Previously, Milne   \cite{m1,m2} had
 obtained the  triangular number identities by a
 third approach.

The present paper rests on two simple observations. First, we
note that the elliptic function identities used
by Zagier can be written as pfaffian evaluations, which follow from
 classical determinant evaluations due to Frobenius and
Stickelberger \cite{fr,fs}. This leads to a new proof of the
Kac--Wakimoto conjecture. Although it is  closely related to
Zagier's proof, it has the advantage of showing how the result could
have been discovered in the nineteenth century, without using
affine superalgebras.

To explain our second observation, we  recall that
the denominator formulas for $Q(2m-1)$ and $Q(2m)$ contain
 $m$ free variables, $x_1,\dots,x_m$.
The Kac--Wakimoto triangular
 number identities are obtained as a limit case when
  $x_j\rightarrow 1$ for all $j$. 
Getz and Mahlburg
\cite{gm} showed that in the case of $Q(2m-1)$,  
letting  $x_j\rightarrow \omega_{2m}^j$, where $\omega_{2m}$ denotes a
primitive $2m$:th 
root of unity, one similarly obtains an identity for
$\triangle_{2m}$, see \eqref{gm} below.  The case $m=1$ gives
back \eqref{t2}, while  $m=2$ gives
\begin{equation}\label{t42}\triangle_{4}(n)=\sum_{\substack{k_1l_1+k_2l_2=8n+4\\ k_i \text{ and } 
l_i \text{ odd positive}\\ k_1\equiv \pm 1,\, k_2\equiv\pm 3\! \mod
  8}}(-1)^{\frac12(k_1-1)+\frac12(k_2-3)}.\end{equation}
In the case of $Q(2m)$,  Getz and Mahlburg  let
  $x_j\rightarrow\omega_{2m+1}^j$. This leads to an
  identity for $\triangle_{2m+1}$, but of a more complicated type
  that we will not consider here.

Our second observation is  that, more
generally, starting from the denominator formula for  $Q(2m-1)$
and letting
$x_j\rightarrow \omega_d^j$, where $d$ is any positive divisor of
$2m$, gives an exact formula for $\triangle_{4m^2/d}$.
 Thus, for any
$m$  we have results corresponding to $d=1$
(Kac--Wakimoto, Milne, Zagier), $d=2$ (new), $d=m$ (new) and $d=2m$ (Getz--Mahlburg),
and if $m$ is neither prime nor equal to $4$ there are additional
results related to the remaining divisors of $2m$.
When $m=2$, we may choose $d=1,2,4$, which apart from \eqref{t16} and
\eqref{t42} gives the, to our knowledge, new identity
$$\triangle_8(n)=\frac{1}{2^5}\left(
\sum_{\substack{k_1l_1+k_2l_2=4n+4\\ k_i \text{ and } l_i \text{
  odd positive}\\ k_1\equiv -k_2\!\mod 4}}(k_1+k_2)^2
-\sum_{\substack{k_1l_1+k_2l_2=4n+4\\ k_i \text{ and } l_i \text{
  odd positive}\\ k_1\equiv k_2\!\mod 4}}(k_1-k_2)^2
\right).
$$

Moreover, in the case of $Q(2m)$, letting
$x_j\rightarrow \omega_d^j$, where $d\mid 2m$ or $d\mid 2m+2$ gives 
an exact formula for $\triangle_{4m(m+1)/d}$, which is of a simpler
type than for the choice $d=2m+1$  in \cite{gm}.
When $m=1$, we may choose $d=1,2,4$, which gives back 
 the three fundamental triangular number
identities \eqref{l}. 
Possibly,
this unified proof of \eqref{l} is  new.  
When $m=2$, the admissible values of $d$ are $1,2,3,4,6$, which apart from \eqref{t24}
 gives the new identities
{\allowdisplaybreaks
\begin{align*}
\triangle_4(n)&=\sum_{\substack{k_1l_1+k_2l_2=6n+3\\
k_i \text{   positive},\  l_i\text{  odd  positive}\\
k_1\equiv\pm 1,\ k_2\equiv\pm 2\!\mod 6
}}(-1)^{\chi(k_1\equiv -1\!\mod 6)+\chi(k_2\equiv -2\!\mod 6)}
,\\
\intertext{where  $\chi(\operatorname{true})=1$, $\chi(\operatorname{false})=0$,}
\triangle_6(n)&=
\frac12\sum_{\substack{k_1l_1+k_2l_2=4n+3\\k_i \text{ positive}, 
\ l_i\text{
    odd  positive}\\ k_1 \text{ odd},\ k_2\equiv 2\!\mod 4}}(-1)^{\frac12(k_1-1)}
k_2,\\
\triangle_8(n)&=\frac1{2\cdot 3^2}\left(
\sum_{\substack{k_1l_1+k_2l_2=3n+3\\k_i\text{ positive},\ l_i\text{
    odd  positive}\\ k_1\equiv -k_2\not\equiv 0\!\mod 3}}(k_1+k_2)^2\right.\\*
&\left.\hspace{5cm}-
\sum_{\substack{k_1l_1+k_2l_2=3n+3\\k_i\text{ positive},\ l_i\text{
    odd  positive}\\ k_1\equiv k_2\not\equiv 0\!\mod3}}(k_1-k_2)^2
\right)
,\\
\triangle_{12}(n)&=\frac 1{2^3}\sum_{\substack{k_1l_1+k_2l_2=2n+3\\
k_1,\,l_1,\,l_2\text{ odd positive}\\ k_2\text{ even positive}}}
k_1k_2^3.\\
\end{align*}
}

The paper is organized as follows. Section \ref{ps} contains
preliminaries on theta functions and pfaffians, in particular our new
proof of the Kac--Wakimoto conjecture. In Section 
\ref{ssc} we review the special case $d=1$ considered by Kac and
Wakimoto. Our main result is given in Theorem \ref{mt},  and in a
slightly different form in Corollary \ref{mtc}; it is proved in
Section \ref{pmts}.

{\bf Acknowledgement:} I thank Eric Rains for illuminating
correspondence on pfaffians.

\section{Preliminaries}
\label{ps}

\subsection{Notation}
The letter $q$ will denote a number such that $0<q<1$, which
will be suppressed from the notation whenever convenient. 
Thus, we write
$$(x)_\infty=(x;q)_\infty=\prod_{j=0}^\infty(1-xq^j).$$
Let $\omega_d=e^{2\pi i/d}$. Then,
\begin{equation}\label{rup}\prod_{k=1}^d(x\omega_d^k;q)_\infty=
(x^d;q^d)_\infty.\end{equation}

We  introduce the theta function
$$\theta(x)=\theta(x;q)=(x,q/x;q)_\infty,$$
which satisfies
\begin{equation}\label{ti}\theta(x^{-1})=\theta(qx)=-x^{-1}\theta(x).
\end{equation}
We will sometimes use the shorthand notation
$$(a_1,\dots,a_n)_\infty=(a_1,\dots,a_n;q)_\infty
=(a_1;q)_\infty\dotsm(a_n;q)_\infty, $$
$$\theta(a_1,\dots,a_n)=\theta(a_1,\dots,a_n;q)
=\theta(a_1;q)\dotsm\theta(a_n;q). $$

We need the classical Laurent expansions
\begin{equation}\label{el}\frac 1x
\frac{\theta(x)}{\theta(\sqrt q x)}=-\frac{(\sqrt
  q)_\infty^2}{(q)_\infty^2} \sum_{k=-\infty}^\infty
  \frac{(\sqrt q x)^k}{1-q^{k+\frac12}},\qquad
  q^\frac12<|x|<q^{-\frac12},
\end{equation}
\begin{equation}\label{ol}x\,\frac{\theta'(x)}{\theta(x)}=-\sum_{k\neq
  0}\frac{x^{k}}{1-q^k},\qquad q<|x|<1,\end{equation}
which can both be derived from Ramanujan's ${}_1\psi_1$ summation \cite{gr}.

For the generating function for   triangular numbers we use
the notation
$$\triangle(q)=\sum_{n=0}^\infty q^{\frac
  12n(n+1)}=\frac12\sum_{n=-\infty}^\infty q^{\frac 12n(n+1)}.  $$
By  Jacobi's triple product identity, we have the product formula
$$
\triangle(q)=(q,-q,-q;q)_\infty=\frac{(q^2;q^2)_\infty}{(q;q^2)_\infty}.
$$
We write $\triangle_k(n)$ for the coefficients in the
Taylor expansion
\begin{equation}\label{tkd}\triangle(q)^k=\sum_{n=0}^\infty
  \triangle_k(n)q^n,\qquad |q|<1,
\end{equation}
which  count
 the  representations of $n$ as a
sum of $k$ triangular numbers. As is customary,
representations obtained from each other by reordering the terms
are considered as different.

\subsection{Pfaffians}

We recall some basic facts about pfaffians.
The pfaffian of a skew-symmetric even-dimensional matrix
$A=(a_{ij})_{i,j=1}^{2m}$ is given by 
$$\pfaff_{1\leq i,j\leq 2m}(a_{ij})
=\frac{1}{2^m m!}\sum_{\sigma\in S_{2m}}\sgn(\sigma)\prod_{i=1}^m
a_{\sigma(2i-1),\sigma(2i)}.$$
Equivalently,
\begin{equation}\label{rpf}\pfaff_{1\leq i,j\leq 2m}(a_{ij})
=\sum_{\sigma\in S_{2m}/G}\sgn(\sigma)\prod_{i=1}^m
a_{\sigma(2i-1),\sigma(2i)},\end{equation}
with $G$ the subgroup of order $2^m m!$ consisting of
 permutations preserving the set of  pairs 
$\{\{1,2\},\{3,4\},\dots,\{2m-1,2m\}\}$. Its main property is
$$\pfaff(A)^2=\det(A).$$

For an
\emph{odd}-dimensional skew-symmetric matrix
$A=(a_{ij})_{i,j=1}^{2m+1}$, we similarly define
$$\pfaff_{1\leq i,j\leq 2m+1}(a_{ij})=\frac{1}{2^m m!}
\sum_{\sigma\in S_{2m+1}}\sgn(\sigma)\prod_{i=1}^m
a_{\sigma(2i-1),\sigma(2i)}.$$
It is easy to check that, with this definition,
$$\pfaff(A)=\pfaff(B), $$ 
where $B=(b_{ij})_{i,j=1}^{2m+2}$ is the matrix
$$B=\left(\begin{matrix}A&\begin{matrix}1\\\vdots \\ 1\end{matrix}
    \\\begin{matrix} -1 & \cdots & -1\end{matrix} & 0\end{matrix}\right). $$

\subsection{The Frobenius determinant and elliptic pfaffians}

Frobenius \cite{fr} obtained the determinant evaluation
\begin{equation}\label{fd}\det_{1\leq i,j\leq
  n}\left(\frac{\theta(tx_iy_j)}{\theta(t,x_iy_j)}\right)
  =\frac{\theta(tx_1\dotsm x_ny_1\dotsm y_n)\prod_{1\leq i<j\leq
  n}x_jy_j\theta(x_i/x_j,y_i/y_j)}{\theta(t)\prod_{i,j=1}^n\theta(x_iy_j)},\end{equation} 
which will form the basis of our analysis. It gives an elliptic
extension of the Cauchy determinant
$$\det_{1\leq i,j\leq
  n}\left(\frac{1}{x_i+y_j}\right)
  =\frac{\prod_{1\leq i<j\leq
  n}(x_i-x_j)(y_i-y_j)}{\prod_{i,j=1}^n(x_i+y_j)}. $$
For other recent applications of \eqref{fd}, see
 \cite{c,kn,r,ru}. The reader interested in elliptic determinant
  evaluations should  consult \cite{kr,rs} for more information 
  and further references.

We also need a  determinant evaluation due to
Frobenius and Stickelberger \cite{fs}, which may be obtained as a
degenerate case of \eqref{fd}. Namely,    rewriting \eqref{fd} as
\begin{multline*}\det_{1\leq i,j\leq
  n+1}\left(\begin{matrix}\displaystyle\frac{\theta(tx_iy_j)}{\theta(t,x_iy_j)}&\begin{matrix}0\\
  \vdots \\ 0\end{matrix}\\\begin{matrix} -1& \cdots & -1\end{matrix} & \theta(t)\end{matrix}\right)\\
  =\frac{\theta(tx_1\dotsm x_ny_1\dotsm y_n)\prod_{1\leq i<j\leq
  n}x_jy_j\theta(x_i/x_j,y_i/y_j)}{\prod_{i,j=1}^n\theta(x_iy_j)}
\end{multline*}
and then subtracting $\theta(t)^{-1}$ times the last row from the
  previous ones gives
\begin{multline*}\det_{1\leq i,j\leq
  n+1}\left(\begin{matrix}\displaystyle\frac{\theta(tx_iy_j)-\theta(x_iy_j)}{\theta(t,x_iy_j)}&\begin{matrix}1\\
  \vdots \\ 1\end{matrix}\\ \begin{matrix} -1 & \cdots & -1 \end{matrix} & \theta(t)\end{matrix}\right)\\
  =\frac{\theta(tx_1\dotsm x_ny_1\dotsm y_n)\prod_{1\leq i<j\leq
  n}x_jy_j\theta(x_i/x_j,y_i/y_j)}{\prod_{i,j=1}^n\theta(x_iy_j)}.
\end{multline*} 
 We  now let $t\rightarrow 1$, obtaining in the limit the
 Frobenius--Stickelberger determinant
\begin{multline}\label{fs}\det_{1\leq i,j\leq
  n+1}\left(\begin{matrix}\displaystyle\frac{x_iy_j\theta'(x_iy_j)}{\theta'(1)\theta(x_iy_j)}&\begin{matrix}1\\
  \vdots \\ 1\end{matrix}\\ \begin{matrix}
-1 & \cdots & -1 \end{matrix} & 0\end{matrix}\right)\\
  =\frac{\theta(x_1\dotsm x_ny_1\dotsm y_n)\prod_{1\leq i<j\leq
  n}x_jy_j\theta(x_i/x_j,y_i/y_j)}{\prod_{i,j=1}^n\theta(x_iy_j)}.
\end{multline}

We are  interested in pfaffian evaluations related to
\eqref{fd} and \eqref{fs}. In \eqref{fd}, we let $n=2m$ and 
choose $t=\sqrt q$,  $y_j=\sqrt q/x_j$. 
Using
\eqref{ti}, the resulting identity can be written as
$$
\det_{1\leq i,j\leq 2m}\left(\frac{\theta(x_j/x_i)}{x_j\theta(\sqrt
    qx_j/x_i)}\right)= q^{\frac12m(m-1)}\prod_{i=1}^{2m}x_i^{2m-2i}\prod_{1\leq i<j\leq
    2m}\frac{\theta(x_j/x_i)^2}{\theta(\sqrt qx_j/x_i)^2}.$$
The matrix on the left is
    skew-symmetric, so we can almost deduce that
\begin{equation}\label{ep}
\pfaff_{1\leq i,j\leq 2m}\left(\frac{\theta(x_j/x_i)}{x_j\theta(\sqrt
    qx_j/x_i)}\right)= q^{\frac
    14m(m-1)}\prod_{i=1}^{2m}x_i^{m-i}\prod_{1\leq i<j\leq 
    2m}\frac{\theta(x_j/x_i)}{\theta(\sqrt qx_j/x_i)}.\end{equation}
More precisely, we know that \eqref{ep} holds up to a factor $\pm 1$
    (possibly depending on $m$). It is not hard
    to show directly that this factor is  always $+1$ (cf.\ the
    final paragraph of \cite{z}), but since that will anyway be clear from our
    computations below,  see Remark \ref{sr}, we will for the moment
    assume that \eqref{ep} is valid.

Applying the same argument to \eqref{fs}, using also 
$\theta'(1)=-(q)_\infty^2$, $\theta(\sqrt q)=(\sqrt q)_\infty^2$,
we obtain, up to a factor $\pm 1$, the
odd-dimensional pfaffian evaluation
\begin{multline}\label{op}
\pfaff_{1\leq i,j\leq 2m+1}\left(\frac{x_i\theta'(\sqrt
  qx_i/x_j)}{x_j\theta(\sqrt q x_i/x_j)}\right)\\
= q^{\frac
  14m(m-1)}\frac{(q)_\infty^{2m}}{(\sqrt
  q)_\infty^{2m}}\prod_{i=1}^{2m+1}x_i^{m+1-i}\prod_{1\leq i<j\leq 
    2m+1}\frac{\theta(x_j/x_i)}{\theta(\sqrt qx_j/x_i)}.\end{multline}
We will see below that the sign chosen in \eqref{op} is correct.

The  evaluations \eqref{ep} and \eqref{op} appear as 
\cite[Eq.\ (7)]{z} (where the pfaffians are written
explicitly as alternating sums). 
As is demonstrated in \cite{z}, they are
equivalent to the denominator formula for the  
super\-algebra $Q(2m-1)$ and $Q(2m)$, respectively,
 which was conjectured by Kac and Wakimoto \cite{kw}.
Thus, our observation that these identities follow from \eqref{fd}
 yields a new proof of the Kac--Wakimoto conjecture.

\begin{remark}
It is interesting to compare \eqref{ep} with some other recent
elliptic pfaffian evaluations.
We first note that letting $n=2m$,
$t=-1$  and
$y_j=-1/x_j$ in \eqref{fd} we obtain, up to a factor $\pm 1$,
\begin{equation}\label{sp}
\pfaff_{1\leq i,j\leq 2m}\left(\frac{\theta(x_j/x_i)}{\theta(-x_j/x_i)}\right)
=\prod_{1\leq i<j\leq 2m}\frac{\theta(x_j/x_i)}{\theta(-x_j/x_i)}.
\end{equation}
For $p=0$, this is  Schur's identity \cite{s}
$$\pfaff_{1\leq i,j\leq 2m}\left(\frac{x_i-x_j}{x_i+x_j}\right)
=\prod_{1\leq i<j\leq 2m}\frac{x_i-x_j}{x_i+x_j};$$
in particular, the sign  in \eqref{sp} is correct.
Alternatively, one may
 obtain \eqref{sp} as the \emph{modular dual} of \eqref{ep}.
Namely, the two results are related by the modular transformation
for Jacobi theta functions, which in our notation takes the form
$$\theta(e^{2\pi ix};\tilde q)=-i\sqrt
h\,\frac{q^{\frac18}(q;q)_\infty}{\tilde q^{\frac18}(\tilde q;\tilde
  q)_\infty}\,e^{\pi ix}q^{\frac12x(x-1)}\,\theta(q^x;q),$$
where
$$q=e^{-2\pi h},\qquad \tilde q=e^{-2\pi/h}.$$
A different  elliptic extension of Schur's identity was recently
obtained by Okada \cite{ok}, namely,
\begin{multline}\label{okp}\pfaff_{1\leq i,j\leq
  2m}\left(\frac{x_i\theta(x_j/x_i,zx_ix_j,wx_ix_j)}
{\theta(x_ix_j,z,w)}\right)\\
=\frac{\theta(zx_1\dotsm x_{2m},wx_1\dotsm
  x_{2m})}{\theta(z,w)}\prod_{1\leq i<j\leq 2m}
\frac{x_i \theta(x_j/x_i)}{\theta(x_ix_j)}. 
\end{multline}
We also mention
Rains'  pfaffian evaluation \cite{r}
  $$\pfaff_{1\leq i,j\leq
    2m}\left(\frac{\theta(x_ix_j,x_j/x_i)}{x_j\theta(\sqrt q x_ix_j,\sqrt
     qx_j/x_i)}\right)=q^{\frac 12m(m-1)}\prod_{1\leq i<j\leq 2m }
\frac{\theta(x_ix_j,x_j/x_i)}{x_j\theta(\sqrt q x_ix_j,\sqrt
     qx_j/x_i)}.$$
Its modular dual is
 $$\pfaff_{1\leq i,j\leq
    2m}\left(\frac{\theta(x_ix_j,x_j/x_i)}{\theta(-
    x_ix_j,-x_j/x_i)}\right)
=\prod_{1\leq i<j\leq 2m }\frac{\theta(x_ix_j,x_j/x_i)}{\theta(-
    x_ix_j,-x_j/x_i)},
$$
which gives a third elliptic extension of Schur's pfaffian, different
    from both \eqref{sp} and \eqref{okp}.
\end{remark}

\begin{remark}
The Frobenius determinant \eqref{fd} has been generalized to higher
genus Riemann surfaces by Fay \cite[Corollary~2.19]{f1}. Similarly, 
it follows from the work of Fay \cite{f2} that the pfaffian evaluations
 \eqref{ep} and
\eqref{sp} can be extended to Prym varieties. We owe this piece of information to Eric Rains.  
\end{remark}

\section{The special case $d=1$}
\label{ssc}

Kac and Wakimoto utilized their (at that time conjectural) denominator
formula to obtain new formulas for the number of representations of an
integer as a sum of $4m^2$ and $4m(m+1)$ triangles.
 Before discussing generalizations, 
 it will be convenient to review the details of this special case.  

We first consider the pfaffian \eqref{ep}. The first step is to use
\eqref{el} to expand the left-hand side as a   multiple Laurent
series. We must then assume
\begin{equation}\label{xc}|x_i/x_j|<q^{-1/2},\qquad i\neq j.\end{equation}
 After interchanging the finite and infinite summations,
the left-hand side of \eqref{ep}  takes the form
\begin{equation}\label{ekw}\frac{(-1)^m(\sqrt q)_\infty^{2m}}{2^m m! (q)_\infty^{2m}}
\sum_{k_1,\dots,k_m=-\infty}^\infty\,\prod_{i=1}^m\frac{q^{\frac 12 k_i}}{1-q^{k_i+\frac
  12}}
\sum_{\sigma\in S_{2m}}\sgn(\sigma)\prod_{i=1}^m
x_{\sigma(2i)}^{k_i}x_{\sigma(2i-1)}^{-k_i-1}.\end{equation}

We next  make the  specialization $x_j=t^j$.
 Then the inner sum in \eqref{ekw}   becomes a
special case of the Vandermonde determinant
\begin{equation}\label{vd}\sum_{\sigma\in S_n}\sgn(\sigma)\prod_{i=1}^n
y_i^{\sigma(i)-1}=\det_{1\leq i,j\leq n}(y_i^{j-1})=\prod_{1\leq
  i<j\leq n}(y_j-y_i).\end{equation}
Namely, choosing $n=2m$,
$(y_1,\dots,y_n)=(t^{-k_1-1},t^{k_1},\dots,t^{-k_m-1},t^{k_m})$, we
obtain after simplification
\begin{multline*}
\sum_{\sigma\in S_{2m}}\sgn(\sigma)\prod_{i=1}^m
t^{k_i\sigma(2i)-(k_i+1)\sigma(2i-1)}\\ 
\begin{split}&=t^{-m}\prod_{i=1}^m(t^{k_i}-t^{-k_i-1})\\
&\quad\times\prod_{1\leq
  i<j\leq
  m}(t^{k_j}-t^{k_i})(t^{k_j}-t^{-k_i-1})(t^{-k_j-1}-t^{k_i})(t^{-k_j-1}-
  t^{-k_i-1})\\
&=(-1)^mt^{-\frac12m(3m+1)}\prod_{i=1}^mt^{(1-2m)k_i}
(1-t^{2k_i+1})\prod_{1\leq
  i<j\leq m} (t^{k_j}-t^{k_i})^2(1-t^{k_i+k_j+1})^2.\end{split}\end{multline*}
This gives
\begin{multline}\label{epx}
\frac{1}{2^m m!}
\sum_{k_1,\dots,k_m=-\infty}^\infty\,\prod_{i=1}^m (t^{1-2m}q^{\frac 12})^{k_i}
\frac{1-t^{2k_i+1}}{1-q^{k_i+\frac
  12}}\prod_{1\leq
  i<j\leq m} (t^{k_j}-t^{k_i})^2(1-t^{k_i+k_j+1})^2\\
= q^{\frac
    14m(m-1)}t^{-\frac 16m(m-1)(4m+1)} 
\frac{(q)_\infty^{2m}}{(\sqrt q)_\infty^{2m}}
\prod_{1\leq i<j\leq 
    2m}\frac{\theta(t^{j-i})}{\theta(\sqrt qt^{j-i})},\end{multline}
which holds for $q^{1/2}<|t^{2m-1}|<q^{-1/2}$. Here we also used
$$t^{\frac12m(3m+1)}\prod_{i=1}^{2m}t^{i(m-i)}=t^{\frac12m(3m+1)-\frac 13m(m+1)(2m+1)}=t^{-\frac 16m(m-1)(4m+1)} .$$

It is clear 
that the left-hand side of \eqref{ekw}, and thus also of \eqref{epx},
is invariant 
under the change of variables $k_i\mapsto -k_i-1$, for any $i$. Thus,
if we multiply  by $2^m$ we may assume that each $k_i$ is
positive. It will also be convenient to make the change of summation
variables  $k_i\mapsto (k_i-1)/2$, giving
\begin{multline}\label{epxo}
\frac{1}{m!}
\sum_{\substack{k_1,\dots,k_m\\ \text{odd positive}}}\,
\prod_{i=1}^m (t^{1-2m}q^{\frac 12})^{\frac12(k_i-1)}
\frac{1-t^{k_i}}{1-q^{\frac 12k_i}}\\
\begin{split}&\quad\times\prod_{1\leq
  i<j\leq m} (t^{\frac12(k_j-1)}-t^{\frac12(k_i-1)})^2(1-t^{\frac12(k_i+k_j)})^2\\
&=  q^{\frac
    14m(m-1)}t^{-\frac 16m(m-1)(4m+1)} 
\frac{(q)_\infty^{2m}}{(\sqrt q)_\infty^{2m}}
\prod_{1\leq i<j\leq 
    2m}\frac{\theta(t^{j-i})}{\theta(\sqrt qt^{j-i})}.
\end{split}\end{multline}

Following Kac and Wakimoto, we now divide both sides of \eqref{epxo}
by
$$\prod_{1\leq i<j\leq 2m}(1-t^{j-i})$$
 and then let $t\rightarrow
1$. On the right-hand side, we have
\begin{equation}\label{rrs}\frac{(q)_\infty^{2m}}{(\sqrt q)_\infty^{2m}}
\prod_{1\leq i<j\leq 
    2m}\frac{\theta(t^{j-i})}{(1-t^{j-i})\theta(\sqrt qt^{j-i})}=
\prod_{i,j=1}^{2m}\frac{(qt^{j-i})_\infty}{(\sqrt
  qt^{j-i})_\infty}\rightarrow \triangle(\sqrt q)^{4m^2}.
 \end{equation}
On the left-hand side, we consider the factor
$$
\frac{\prod_{i=1}^m (1-t^{k_i})\prod_{1\leq
  i<j\leq m} (t^{\frac12(k_j-1)}-t^{\frac12(k_i-1)})^2
(1-t^{\frac12(k_i+k_j)})^2}{\prod_{1\leq i<j\leq 2m}(1-t^{j-i})}.
$$
Since $m+4\binom m2=\binom{2m}2$, the denominator and numerator
vanish of the same order at $t=1$, so that the quotient tends to
$$
\frac{\prod_{i=1}^m k_i\prod_{1\leq
  i<j\leq m} (\frac{k_j-k_i}{2})^2
(\frac{k_i+k_j}{2})^2}{\prod_{1\leq i<j\leq 2m}(j-i)}
=\frac{\prod_{i=1}^m k_i\prod_{1\leq
  i<j\leq m} (k_j^2-k_i^2)^2}{4^{m(m-1)}\prod_{j=1}^{2m-1}j!}.
$$
Thus, we obtain the identity
$$\frac{q^{-\frac14m(m-1)}}{4^{m(m-1)}m!\prod_{j=1}^{2m-1}j!}
\sum_{\substack{k_1,\dots,k_m\\ \text{odd positive}}}
\,\prod_{i=1}^m 
\frac{q^{\frac14(k_i-1)}k_i}{1-q^{\frac 12k_i}}\prod_{1\leq
  i<j\leq m} (k_j^2-k_i^2)^2
=\triangle(\sqrt q)^{4m^2}.$$

We now expand the left-hand side directly as a power series in
$\sqrt q$, using
\begin{equation}\label{gs}\frac 1{1-q^{\frac 12k_i}}=\sum_{l_i \text{ odd positive}}q^{\frac
  14k_i(l_i-1)},\end{equation}
which gives
$$
\frac{1}{4^{m(m-1)}m!\prod_{j=1}^{2m-1}j!}
\sum_{\substack{k_1,\dots,k_m, l_1,\dots,l_m\\ \text{odd positive}}}
q^{\frac 14(k_1l_1+\dotsm+k_ml_m-m^2)}\prod_{i=1}^m k_i\prod_{1\leq
  i<j\leq m} (k_j^2-k_i^2)^2.
$$
In conclusion, this proves that
\begin{equation}\label{kmt}\triangle_{4m^2}(n)=\frac{1}{4^{m(m-1)}m!\prod_{j=1}^{2m-1}j!}
\sum_{\substack{k_1l_1+\dotsm+k_ml_m=2n+m^2\\ k_i \text{ and } l_i \text{ odd positive}}}
\,\prod_{i=1}^m k_i\prod_{1\leq
  i<j\leq m} (k_j^2-k_i^2)^2.\end{equation}

\begin{remark}\label{sr}
Had we started from the identity similar to \eqref{ep}
but with the right-hand side multiplied by $-1$, we would have
obtained a similarly modified version of \eqref{kmt}, which would
clearly be absurd. Thus, the sign chosen in \eqref{ep} is correct.  
\end{remark}

We now turn to the case of  \eqref{op}, where we will be less
detailed. However, we write down some intermediate steps for later
reference.  Still assuming \eqref{xc}, we  apply
 \eqref{ol} to rewrite the left-hand-side of \eqref{op} as
\begin{equation}\label{okw}\frac{(-1)^m q^{-\frac 12m}}{2^m m!}
\sum_{k_1,\dots,k_m\neq 0}\,\prod_{i=1}^m\frac{q^{k_i/2}}{1-q^{k_i}}
\sum_{\sigma\in S_{2m+1}}\sgn(\sigma)\prod_{i=1}^m
x_{\sigma(2i-1)}^{k_i}x_{\sigma(2i)}^{-k_i}.\end{equation}

As before, we  choose $x_i=t^i$.
By the case $n=2m+1$,
$(y_1,\dots,y_{n})=(t^{k_1},t^{-k_1},\dots, t^{k_m},t^{-k_m},1)$ of
\eqref{vd}, the inner sum in \eqref{okw} equals
$$(-1)^m\prod_{i=1}^mt^{-2mk_i}(1-t^{k_i})^2(1-t^{2k_i})\prod_{1\leq i<j\leq
m}(t^{k_j}-t^{k_i})^2(1-t^{k_i+k_j})^2. $$
Exploiting the symmetry $k_i\mapsto -k_i$, we  reduce the summation
to positive $k_i$, giving 
\begin{multline}\label{opx}
\frac{1}{m!}\sum_{k_1,\dots,k_m=1}^\infty \,
\prod_{i=1}^m\frac{(q^{\frac12}t^{-2m})^{k_i}}{1-q^{k_i}}
(1-t^{k_i})^2(1-t^{2k_i})\\
\begin{split}&\quad\times\prod_{1\leq i<j\leq
m}(t^{k_j}-t^{k_i})^2(1-t^{k_i+k_j})^2\\
&= q^{\frac
  14m(m+1)}t^{-\frac13m(m+1)(2m+1)}\frac{(q)_\infty^{2m}}{(\sqrt
  q)_\infty^{2m}}\prod_{1\leq i<j\leq 
    2m+1}\frac{\theta(t^{j-i})}{\theta(\sqrt q t^{j-i})},
\end{split}\end{multline}
which holds for $q^{\frac 12}<|t^{2m}|<q^{-\frac 12}$.

Dividing \eqref{opx} by $\prod_{1\leq i<j\leq   2m+1}(1-t^{j-i})$
and letting $t$ tend to $1$ gives
$$\frac{ 2^mq^{-\frac14m(m+1)}}{m!\prod_{j=1}^{2m}j!}
\sum_{k_1,\dots,k_m=1}^\infty \,
\prod_{i=1}^m 
\frac{q^{\frac12k_i}k_i^3}{1-q^{k_i}}\prod_{1\leq
  i<j\leq m} (k_j^2-k_i^2)^2
= \triangle(\sqrt q)^{4m(m+1)}.$$
Expanding the denominator using
\begin{equation}\label{ope}\frac 1{1-q^{k_i}}=\sum_{l_i \text{\ odd positive}}q^{\frac12k_i(l_i-1)}\end{equation} 
and identifying the coefficient of $q^{n/2}$ we obtain 
\begin{equation}\label{kmt2}\triangle_{4m(m+1)}(n)=
\frac{ 2^m}{m!\prod_{j=1}^{2m}j!}\sum_{\substack{k_1l_1+\dotsm
  +k_ml_m=n+\frac 12m(m+1)\\ k_i\text{ positive},\ l_i \text{
  odd  positive}}}\,\prod_{i=1}^m k_i^3\prod_{1\leq
  i<j\leq m} (k_j^2-k_i^2)^2.\end{equation}
In particular, we conclude that the choice of sign in \eqref{op} is correct.

\section{The general case}
\label{sd}

 When $(k_1,\dots,k_m)$ and $(l_1,\dots,l_m)$  are
multi-indices, let us write  
$$(k_1,\dots,k_m)\simeq(l_1,\dots,l_m) \mod n$$
if they are equal modulo $n$ up to reordering and sign, 
that is, if there exists a permutation 
 $\sigma\in S_m$ and numbers $\varepsilon_i\in\{\pm 1\}$ such that 
 $k_{\sigma(i)}\equiv\varepsilon_i l_i\pod n$ for $i=1,\dots, m$. 
In this notation, we can state our main result as follows. 

\begin{theorem}
\label{mt}
Let $d$, $m$ and $x$ be non-negative integers.
Then, if $d\mid 2m$,
\begin{multline}\label{mt1}\sum_{\substack{k_1l_1+\dots+k_ml_m=m^2+x\\
    k_i \text{ \emph{and} } l_i \text { \emph{odd positive}}\\ 
(k_1,\dots,k_m)\simeq(1,3,5,\dots,2m-1)\ (2d)}}
(-1)^{|\{i;\,d+1\leq k_i\leq
  2d-1\!\mod 2d\}|}\\
\begin{split}&\quad\times\prod_{\substack{1\leq i\leq m\\
k_i\equiv d\ (2d)}} k_i\prod_{\substack{1\leq
  i<j\leq m\\k_i\equiv k_j\,(2d)}}{\left(\frac{k_j-k_i}{2}\right)^2}
\prod_{\substack{1\leq
  i<j\leq m\\k_i\equiv -k_j\,(2d)}}
{\left(\frac{k_j+k_i}{2}\right)^2}\\
&= (-1)^{\binom d2\binom{2m/d}2}d^{(2m-d)m/d}m!
\prod_{l=1}^{(2m-d)/d}l!^d
\,
\triangle_{4m^2/d}(x/2d).
\end{split} \end{multline}
Moreover, if $d\mid 2m$ or $d\mid 2m+2$,
\begin{multline}\label{mt2}
\sum_{\substack{k_1l_1+\dotsm+ k_ml_m=\frac12m(m+1)+x\\k_i \text{
      \emph{positive}},\ l_i \text{ \emph{odd
      positive}}\\(k_1,\dots,k_m)\simeq(1,2,\dots,m)\,(d)}}\,
(-1)^{|\{i;\,(d+1)/2\leq k_i\leq
  d-1\!\mod d\}|}
\\
\begin{split}&\quad\times\prod_{\substack{1\leq i\leq m\\k_i\equiv 0\,(d)}}2k_i^3
\prod_{\substack{1\leq i\leq m\\k_i\equiv d/2\,(d)}}2k_i
\prod_{\substack{1\leq
  i<j\leq m\\k_i\equiv k_j\,(d)}}
(k_j-k_i)^2
\prod_{\substack{1\leq
  i<j\leq m\\k_i\equiv -k_j\,(d)}}
(k_i+k_j)^2,
\\
&=C\,d^{(2m+2-d)m/d}
m!\,
\triangle_{4m(m+1)/d}(x/d),\end{split}\end{multline}
where
\begin{equation}\label{c}C=\begin{cases}(-1)^{\binom{d-1}2\binom{2m/d}{2}} (2m/d)!
\prod_{l=1}^{(2m-d)/d}l!^d, &d\mid 2m,\\
(-1)^{\binom{d-1}2\binom{2(m+1)/d}{2}}
((2m+2-d)/d)!^{-1}
\prod_{l=1}^{(2m+2-d)/d}l!^d, & d\mid 2m+2. 
 \end{cases} \end{equation}
The right-hand side of \eqref{mt1} and \eqref{mt2} should be
interpreted as zero if $2d\nmid x$ and $d\nmid x$, respectively.
\end{theorem}

\begin{remark}\label{r}
Since the sums in Theorem \ref{mt} are symmetric in $k_i$ and vanish
if $k_i=k_j$ for $i\neq j$, each term is repeated $m!$ times. If one
wants to compute the sums, one should first get rid of this
redundancy. This can be done, for instance, by imposing the condition
$k_1<k_2<\dots<k_m$ and deleting the factor $m!$ from the right-hand
side. However, it is  more convenient to, in the case of \eqref{mt1},
first impose the condition  $k_i\equiv\pm(2i-1)\pod{2d}$  
and then the condition $k_i<k_j$
if $i<j$ and $k_i\equiv\pm k_j\pod{2d}$, 
and similarly for \eqref{mt2}. 
\end{remark}

In  Theorem \ref{mt} we have tried to state the results in a unified
form. However, this hides some  structural differences between 
 even and odd $d$. 
In  Corollary \ref{mtc} we rewrite the identities in a 
way that emphasizes these differences, and is also better suited for
computation as indicated in Remark \ref{r}. Moreover, we have removed a
power of $2$ from the left-hand sides, using \eqref{2a} and \eqref{2b} below.

\begin{corollary}
\label{mtc}
Let $d$, $m$ and $x$ be non-negative integers.  Then, if $d$ is odd
and $d\mid m$,
\begin{subequations}\label{mt1ab}
\begin{multline}\label{mt1a}\sum_{\substack{k_1l_1+\dots+k_ml_m=m^2+x\\
    k_i \text{ \emph{and} } l_i \text { \emph{odd positive}}\\
    k_i\equiv \pm(2i-1)\!\mod 2d\\ k_i<k_j \text{ \emph{if} } i<j
    \text{ \emph{and} } k_i\equiv\pm k_j\,(2d)}}
(-1)^{|\{i;\,k_i\equiv
  d+2,d+4,\dots,2d-1\,(2d)\}|}\\
\begin{split}&\quad\times\prod_{\substack{1\leq i\leq m\\
k_i\equiv d\,(2d)}} k_i\prod_{\substack{1\leq
  i<j\leq m\\k_i\equiv k_j\,(2d)}}{(k_j-k_i)^2}
\prod_{\substack{1\leq
  i<j\leq m\\k_i\equiv -k_j\,(2d)}}
{(k_j+k_i)^2}\\
&=(-1)^{\frac {d-1}2\frac md}(2^{2m-d-1}d^{2m-d})^{m/d}
\prod_{l=1}^{(2m-d)/d}l!^d\,\triangle_{4m^2/d}(x/2d),
\end{split} \end{multline}
while if $d$ is  even and $d\mid 2m$,  
\begin{multline}\label{mt1b}\sum_{\substack{k_1l_1+\dots+k_ml_m=m^2+x\\
    k_i \text{ \emph{and} } l_i \text { \emph{odd positive}}\\
    k_i\equiv\pm(2i-1)\!\mod 2d\\ k_i<k_j \text{ \emph{if} } i<j
    \text{ \emph{and} } k_i\equiv\pm k_j\,(2d)}}
(-1)^{|\{i;\,k_i\equiv
  d+1,d+3,\dots,2d-1\,(2d)\}|}\\
\begin{split}&\quad\times\prod_{\substack{1\leq
  i<j\leq m\\k_i\equiv k_j\,(2d)}}{(k_j-k_i)^2}
\prod_{\substack{1\leq
  i<j\leq m\\k_i\equiv -k_j\,(2d)}}
{(k_j+k_i)^2}\\
&=(-1)^{\frac d2\binom{2m/d}{2}}(2d)^{(2m-d)m/d} \prod_{l=1}^{(2m-d)/d}l!^d\,\triangle_{4m^2/d}(x/2d).
\end{split} \end{multline}
\end{subequations}
Moreover, if $d$ is odd and $d\mid m$ or $d\mid m+1$,
\begin{subequations}\label{mt2ab}
\begin{multline}\label{mt2a}
\sum_{\substack{k_1l_1+\dotsm+ k_ml_m=\frac12m(m+1)+x\\k_i \text{
      \emph{positive}},\ l_i \text{ \emph{odd
      positive}}\\k_i\equiv\pm i\!\mod d\\ 
k_i<k_j \text{ \emph{if} } i<j
    \text{ \emph{and} } k_i\equiv\pm k_j\,(d)}}\,
(-1)^{|\{i;\,k_i\equiv (d+1)/2,(d+3)/2,\dots,d-1\,(d)\}|}
\\
\begin{split}&\quad\times
\prod_{\substack{1\leq i\leq m\\k_i\equiv 0\,(d)}}k_i^3
\prod_{\substack{1\leq
  i<j\leq m\\k_i\equiv k_j\,(d)}}
(k_j-k_i)^2
\prod_{\substack{1\leq
  i<j\leq m\\k_i\equiv -k_j\,(d)}}
(k_i+k_j)^2,
\\
&=A\,d^{(2m+2-d)m/d}
\triangle_{4m(m+1)/d}(x/d),\end{split}\end{multline}
where
$$A=\begin{cases}(-1)^{\binom{d-1}2\binom{2m/d}{2}} 
2^{-m/d}(2m/d)!
\prod_{l=1}^{(2m-d)/d}l!^d, &d\mid m,\\
(-1)^{\binom{d-1}2\binom{2(m+1)/d}{2}}
2^{-(m+1-d)/d}&\\
\qquad\times ((2m+2-d)/d)!^{-1}
\prod_{l=1}^{(2m+2-d)/d}l!^d, & d\mid m+1,\end{cases} $$
while if $d$ is even and $d\mid 2m$ or $d\mid 2m+2$, 
\begin{multline}\label{mt2b}
\sum_{\substack{k_1l_1+\dotsm+ k_ml_m=\frac12m(m+1)+x\\k_i \text{
      \emph{positive}},\ l_i \text{ \emph{odd
      positive}}\\k_i\equiv\pm i\!\mod d\\ 
k_i<k_j \text{ \emph{if} } i<j
    \text{ \emph{and} } k_i\equiv\pm k_j\,(d)}}\,
(-1)^{|\{i;\,k_i\equiv(d+2)/2,(d+4)/2,\dots,d-1\,(d)\}|}
\\
\begin{split}&\quad\times
\prod_{\substack{1\leq i\leq m\\k_i\equiv 0\,(d)}}k_i^3
\prod_{\substack{1\leq i\leq m\\k_i\equiv d/2\,(d)}}k_i
\prod_{\substack{1\leq
  i<j\leq m\\k_i\equiv k_j\,(d)}}
(k_j-k_i)^2
\prod_{\substack{1\leq
  i<j\leq m\\k_i\equiv -k_j\,(d)}}
(k_i+k_j)^2,
\\
&=B\,d^{(2m+2-d)m/d}
\triangle_{4m(m+1)/d}(x/d),\end{split}\end{multline}
where
$$B=\begin{cases}(-1)^{\binom{d-1}2\binom{2m/d}{2}} 
2^{-2m/d}(2m/d)!
\prod_{l=1}^{(2m-d)/d}l!^d, &d\mid 2m,\\
(-1)^{\binom{d-1}2\binom{2(m+1)/d}{2}}
2^{-(2m+2-d)/d}&\\
\qquad\times((2m+2-d)/d)!^{-1}
\prod_{l=1}^{(2m+2-d)/d}l!^d, & d\mid 2m+2.\end{cases} $$
\end{subequations}
\end{corollary}

Apart from the case $d=1$ discussed in Section \ref{ssc}, the other extremal
case, $d=2m$ of \eqref{mt1} (or \eqref{mt1b})
and $d=2m+2$ of \eqref{mt2} (or \eqref{mt2b}) is
of special interest. Then, all products on the left are
empty except for the power of
$-1$, so that both sides of the identity have a clear combinatorial
meaning.  
In the case of \eqref{mt1b}, we recover
\begin{equation}\label{gm}\sum_{\substack{k_1l_1+\dots+k_ml_m=m^2+x\\
    k_i \text{ {and} } l_i \text { {odd positive}}\\
    k_i\equiv\pm(2i-1)\!\mod 4m}}
(-1)^{|\{i;\,k_i\equiv
  1-2i\ (4m)\}|}=\triangle_{2m}(x/4m),
 \end{equation}
 which is equivalent to \cite[Corollary 1.3]{gm}.
Similarly, \eqref{mt2b} gives
\begin{equation}\label{gmo}\sum_{\substack{k_1l_1+\dotsm+ k_ml_m=\frac12m(m+1)+x\\k_i \text{
      {positive}},\ l_i \text{ {odd
      positive}}\\k_i\equiv\pm i\!\mod 2m+2}}\,
(-1)^{|\{i;\,k_i=-i\ (2m+2)\}|}=\triangle_{2m}(x/(2m+2)).\end{equation}
The fact  
that these sums vanish for  $4m\nmid x$ and $2m+2\nmid x$,
      respectively, and are 
 otherwise non-negative 
 gives some non-trivial information, see  \cite{gm} for the case of \eqref{gm}.

\section{Proof of Theorem \ref{mt}}
\label{pmts}

To prove Theorem \ref{mt},  we will let
$t\rightarrow\omega_d=e^{2\pi i/d}$ in \eqref{epxo} and \eqref{opx},
 after dividing both sides by a suitable factor. 
We will  assume $d\mid 2m$ in the case of \eqref{epxo} and $d\mid
2m$ or  $d\mid 2m+2$ in the case of 
\eqref{opx}, obtaining \eqref{mt1} and \eqref{mt2},
respectively. As is already clear from \cite{gm},  other
values of $d$ may also have arithmetic consequences, but we restrict here
to the  simplest situation.

Before working out the details, we collect some elementary but useful facts. Note that part (b) of the following Lemma gives a more
explicit description of the range of summation in \eqref{mt1}.

{\allowdisplaybreaks

\begin{lemma}\label{l1}
Assume that $d\mid 2m$ and consider the sequence 
$(1,3,5,\dots,2m-1)$ reduced modulo $2d$. Then:
\renewcommand{\labelenumi}{(\alph{enumi})}
\begin{enumerate}
\item The number of elements of the sequence that are congruent to $i$
  modulo $2d$ equals
\begin{align*}&\frac md+\frac 12, & i&=1,3,\dots,d-1,\\
&\frac md-\frac 12, & i&=d+1,d+3,\dots,2d-1
\end{align*} 
if $2m/d$ is odd (and thus $d$ is even) and 
 $$m/d, \qquad i=1,3,\dots,2d-1 $$
if $2m/d$ is even.
\item
The number of elements of the sequence congruent to $\pm i$
  modulo $2d$ equals 
\begin{align*}&2m/d, & i&=1,3,\dots,2d-1,\ i\neq d,\\
&m/d, & i&=d,\ d \text{ \emph{odd}}.
\end{align*} 
\item
 The number of elements $x$ of the sequence  
 such that $d+1\leq x\leq 2d-1$ modulo $2d$ has the
same parity as $\binom{d}{2}\binom{2m/d}{2}$.
\end{enumerate}
\end{lemma}

\begin{proof}
The proof of  (a) is trivial, and   (b) follows
immediately from (a). As for the final statement, it
follows from   (a) that
the number of such $x$ equals
$$\frac d2\left(\frac md-\frac 12\right),\quad
\frac d2\cdot\frac md , \quad \frac{d-1}{2}\cdot\frac{m}{d}
 $$
according to whether $2m/d$ is odd and $d$ even, $2m/d$  and $d$ are
both even or $2m/d$ is even and $d$ odd, respectively. The quotients of 
 $\binom{d}{2}\binom{2m/d}{2}$ by these numbers equal
$$\frac{2m}{d}(d-1),\quad \left(\frac{2m}{d}-1\right)(d-1),\quad
\left(\frac{2m}{d}-1\right)d, $$
which in each case is  odd.
\end{proof}

In the case of \eqref{mt2}, the corresponding facts are  somewhat more
tedious  to state.

\begin{lemma}\label{l2}
Assume that $d\mid 2m$ or $d\mid 2m+2$ and consider the sequence 
$(1,2,3,\dots,m)$ reduced modulo $d$. Then:
\renewcommand{\labelenumi}{(\alph{enumi})}
\begin{enumerate}
\item The number of elements of the sequence congruent to $i$
  modulo $d$ equals
\begin{align*}&\frac md+\frac 12, & i&=1,2,\dots,d/2,\\
&\frac md-\frac 12, & i&=\frac d2+1,\frac d2+2,\dots,d
\end{align*}
if $2m/d$ is odd; 
 $$m/d, \qquad i=1,2,\dots,d $$
if $2m/d$ is even;
\begin{align*}&\frac {m+1}d+\frac 12, & i&=1,2,\dots,\frac d2-1,\\
 &\frac{m+1}d-\frac 12, & i&=\frac d2,\frac d2+1,\dots,d
\end{align*}
if $(2m+2)/d$ is odd and
\begin{align*}&(m+1)/d, & i&=1,2,\dots,d-1,\\
 &\frac{m+1}d-1, & i&=d
\end{align*}
if $(2m+2)/d$ is even.
\item
The number of elements of the sequence congruent to $\pm i$
  modulo $ d$ equals
\begin{align*}
&2m/d, & i&=1,2,\dots, d-1,\, i\neq d/2,\\
&\frac md+\frac12, & i&=d/2,\\
&\frac md-\frac12, & i&=d
\end{align*}
if $2m/d$ is odd;
\begin{align*}
&2m/d, & i&=1,2,\dots, d-1,\, i\neq d/2,\\
&m/d, & i&=d/2\ (d \text{\emph{ even}}),\, i=d
\end{align*}
if $2m/d$ is even;
\begin{align*}
&2(m+1)/d, & i&=1,2,\dots, d-1,\, i\neq d/2,\\
&\frac{m+1}d-\frac12, & i&=d/2,\,i=d
\end{align*}
if $(2m+2)/d$ is odd and
\begin{align*}
&2(m+1)/d, & i&=1,2,\dots, d-1,\, i\neq d/2,\\
&(m+1)/d, & i&=d/2\ (d \text{\emph{ even}}),\\
&\frac{m+1}d-1, & i&=d
\end{align*}
if $(2m+2)/d$ is even.
\item
 The number of elements $x$  of the sequence
 such that $(d+1)/2\leq x\leq d-1$ modulo $d$ has the
same parity as $\binom{d-1}{2}\binom{2m/d}{2}$ if $d\mid 2m$ and
$\binom{d-1}{2}\binom{2(m+1)/d}{2}$ if  $d\mid 2m+2$.
\end{enumerate}
\end{lemma}

}

{\bf Proof of \eqref{mt1}.}
We first assume $d\mid 2m$ and
consider  the right-hand side of \eqref{epxo}. After dividing
by the factor  
$$P(t)=\prod_{1\leq i<j\leq 2m}(1-t^{j-i}),$$
 we have as in \eqref{rrs}
\begin{multline}\label{rhl}\frac{(q)_\infty^{2m}}{(\sqrt q)_\infty^{2m}}
\prod_{1\leq i<j\leq 
    2m}\frac{\theta(t^{j-i})}{(1-t^{j-i})\theta(\sqrt qt^{j-i})}\Bigg|_{t=\omega_d}=
\prod_{i,j=1}^{2m}\frac{(q\omega_d^{j-i})_\infty}{(\sqrt
  q\omega_d^{j-i})_\infty}
\\
=\left(\prod_{i,j=1}^{d}\frac{(q\omega_d^{j-i})_\infty}{(\sqrt
  q\omega_d^{j-i})_\infty}\right)^{4m^2/d^2}=\left(\prod_{k=1}^{d}\frac{(q\omega_d^{k})_\infty}{(\sqrt
  q\omega_d^{k})_\infty}\right)^{4m^2/d}=\triangle(q^{d/2})^{4m^2/d},
 \end{multline}
where we used \eqref{rup} 
in the last step.

We  compute the   multiplicity of $t=\omega_d$ as a zero
of the right-hand side of \eqref{epxo}, or equivalently of $P$. 
This is the  number of pairs $(i,j)$
such that $1\leq i<j\leq 2m$ and $i\equiv j\pod d$. If we assume
  $j=i+ld$,  $1\leq l\leq (2m-d)/d$, there are $2m-ld$ such pairs. Summing over
 $l$ gives the total multiplicity
\begin{equation}\label{pzm}(2m-d)+(2m-2d)+\dots+2d+d=\frac{m(2m-d)}{d}.
\end{equation}

Next we turn to the 
left-hand side of \eqref{epxo}.

\begin{lemma}\label{l3}
If $d\mid 2m$ and $k_1,\dots, k_m$ are odd, then the multiplicity of 
 $t=\omega_d$ as a zero of
$$\prod_{i=1}^m(1-t^{k_i})\prod_{1\leq
  i<j\leq m} (t^{\frac12(k_j-1)}-t^{\frac12(k_i-1)})^2
(1-t^{\frac12(k_i+k_j)})^2$$
is at least $m(2m-d)/d$, with equality if and
only if 
\begin{equation}\label{kc}(k_1,\dots,k_m)\simeq (1,3,\dots,2m-1)\mod
  2d.\end{equation} 
\end{lemma}

\begin{proof}
Let $a_i$ be the number of $j$, $1\leq
j\leq m$, such that $k_j\equiv\pm i\pod{2d}$. Here we take $i=1,3,5,\dots,d$
  or $i=1,3,5,\dots,d-1$ according to whether $d$ is odd or even.
We note that there are  single zeroes when 
  $k_i\equiv d\pod{2d}$ and double zeroes when $k_i\equiv \pm
  k_j\pod{2d}$, $i<j$.
Thus, if $d$ is odd,
the total multiplicity is
$$a_{d}+2\binom{a_1}{2}+2\binom{a_3}{2}+\dots+2\binom{a_{d-2}}{2}+
4\binom{a_d}{2}, $$
while if $d$ is even, it is
$$2\binom{a_1}{2}+2\binom{a_3}{2}+\dots+2\binom{a_{d-1}}{2}. $$
We want to  minimize these expressions under the condition $\sum_i a_i=m$.
Using, for instance,   Lagrange multipliers, one
checks that in both cases the minimum is $m(2m-d)/d$.
Moreover, it is achieved precisely if
$$ a_1=a_3=\dots=a_{d-2}=2m/d,\qquad
a_{d}=m/d,
$$
for odd $d$ and if
$$ a_1=a_3=\dots=a_{d-1}=2m/d$$
for even $d$. By 
Lemma \ref{l1}.b,  this is in both cases equivalent to \eqref{kc}.
\end{proof}

We may now let $t\rightarrow\omega_d$ termwise in  \eqref{epxo},
concluding that
$$
\frac{1}{m!}\,\sideset{}{'}
\sum_{k_1,\dots,k_m}\,
\lim_{t\rightarrow\omega_d}\frac{T_{k_1,\dots,k_m}(t)}{P(t)}
=  q^{\frac
    14m(m-1)}\omega_d^{-\frac16m(m-1)(4m+1)}
\triangle(q^{d/2})^{4m^2/d}
,$$
where $\sum'$ indicates that the summation variables are positive odd
integers  satisfying \eqref{kc}, and where
\begin{multline*}\begin{split}T_{k_1,\dots,k_m}(t)
&=
\prod_{i=1}^m (t^{1-2m} q^{\frac 12})^{\frac12(k_i-1)}
\frac{1-t^{k_i}}{1-q^{\frac 12k_i}}\\
&\quad\times\prod_{1\leq
  i<j\leq m}
(t^{\frac12(k_j-1)}-t^{\frac12(k_i-1)})^2(1-t^{\frac12(k_i+k_j)})^2.
\end{split}\end{multline*}

We factor $T=T_{k_1,\dots,k_m}(t)$ as $T^1T^2$, where
\begin{align*}
T^1&=\prod_{i=1}^m (t^{1-2m})^{\frac12(k_i-1)}
\prod_{\substack{1\leq i<j\leq
    m\\k_i\equiv k_j\,(2d)}}t^{k_i-1}
\prod_{\substack{1\leq i\leq m\\k_i\not\equiv d\,(2d)}}
(1-t^{k_i})\\
&\quad\times\prod_{\substack{1\leq
  i<j\leq m\\k_i\not\equiv k_j\,(2d)}}
(t^{\frac12(k_i-1)}-t^{\frac12(k_j-1)})^2
\prod_{\substack{1\leq
  i<j\leq m\\k_i\not\equiv -k_j\,(2d)}}
(1-t^{\frac12(k_i+k_j)})^2,\\
T^2&= \prod_{i=1}^m  \frac{q^{\frac14(k_i-1)}}{1-q^{\frac 12k_i}}
\prod_{\substack{1\leq i\leq m\\k_i\equiv d\,(2d)}}
(1-t^{k_i})\\*
&\quad\times\prod_{\substack{1\leq
  i<j\leq m\\k_i\equiv k_j\,(2d)}}
(1-t^{\frac12(k_j-k_i)})^2
\prod_{\substack{1\leq
  i<j\leq m\\k_i\equiv -k_j\,(2d)}}
(1-t^{\frac12(k_i+k_j)})^2.
\end{align*}
Similarly, we write $P=P^1P^2$, where
{\allowdisplaybreaks
\begin{align*}
P^1(t)&=\prod_{\substack{1\leq i<j\leq 2m\\i\not\equiv
    j\,(d)}}(1-t^{j-i}),\\
P^2(t)&=\prod_{\substack{1\leq i<j\leq 2m\\i\equiv
    j\,(d)}}(1-t^{j-i}).\end{align*}
}
Then,
\begin{multline}\label{epxo2}
\frac{1}{m!}\,\sideset{}{'}
\sum_{k_1,\dots,k_m}\,\frac{T^1_{k_1,\dots,k_m}(\omega_d)}{P^1(\omega_d)}
\lim_{t\rightarrow\omega_d}\frac{T^2_{k_1,\dots,k_m}(t)}{P^2(t)}\\
=  q^{\frac
    14m(m-1)}\omega_d^{-\frac16m(m-1)(4m+1)}\triangle(q^{d/2})^{4m^2/d}
,\end{multline}
where $T^1/P^1$ simplifies in view of the following lemma.

\begin{lemma}\label{tpl}
In the notation above,
\begin{equation}\label{tp1}\frac{T^1_{k_1,\dots,k_m}(\omega_d)}{P^1(\omega_d)}=(-1)^{|\{i;\,
    d+1\leq k_i\leq 2d-1\, (2d)\}|
+\binom{d}{2}\binom{2m/d}{2}}\omega_d^{-\frac 16m(m-1)(4m+1)}.
\end{equation}
\end{lemma}

\begin{proof}
Let 
$\tau(k_1,\dots,k_m)
=T^1_{k_1,\dots,k_m}(\omega_d)$. Then, $\tau$ is
visibly symmetric in the $k_i$ and invariant under
 $k_i\mapsto k_i+2d$ for any $i$. It is also easy
to check that
$$\frac{\tau(-k_1,k_2,\dots,k_m)}{\tau(k_1,k_2,\dots,k_m)}= -1, \qquad
k_1\not\equiv 0\pod d. 
 $$
Thus,  $\tau(k_1,\dots,k_m)$ equals, up
to a factor $\pm 1$, $\tau(1,3,5,\dots,2m-1)$. The sign may be computed
using  Lemma \ref{l1}.c, giving
$$\tau(k_1,\dots,k_m)=(-1)^{|\{i;\,d+1\leq k_i\leq 2d-1\, (2d)\}|
+\binom{d}{2}\binom{2m/d}{2}}
\tau(1,3,5,\dots,2m-1).$$

Next, we write
$$P^1(t)=\prod_{\substack{1\leq i<j\leq m\\i\not\equiv
    j\,(d)}}(1-t^{j-i})\prod_{\substack{m+1\leq i<j\leq 2m\\i\not\equiv
    j\,(d)}}(1-t^{j-i})\prod_{\substack{1\leq i\leq m,\,m+1\leq j\leq 2m\\i\not\equiv
    j\,(d)}}(1-t^{j-i}).$$
In the second  product, we replace $(i,j)\mapsto(m+i,m+j)$ and in
    the third  product $(i,j)\mapsto(m+1-i,m+j)$, giving
\begin{equation*}\begin{split}P^1(t)&
=\prod_{\substack{1\leq i<j\leq m\\i\not\equiv
    j\,(d)}}(1-t^{j-i})^2
\prod_{\substack{1\leq i\leq m,\,1\leq j\leq m\\i+j\not\equiv
    1\,(d)}}(1-t^{i+j-1})\\
&=\prod_{\substack{1\leq i<j\leq m\\i\not\equiv
    j\,(d)}}(1-t^{j-i})^2
\prod_{\substack{1\leq i<j\leq m\\i+j\not\equiv
    1\,(d)}}(1-t^{i+j-1})^2\prod_{\substack{1\leq i\leq
    m\\2i\not\equiv 1\, (d)}}(1-t^{2i-1})^2.
\end{split}\end{equation*}
Comparing with the definition of $T^1$ we find that, in general,
$$\frac{T^1_{1,3,\dots,2m-1}(t)}{P^1(t)}=\prod_{i=1}^mt^{(1-2m)(i-1)} \prod_{1\leq i<j\leq m}t^{2i-2}=t^{-\frac16m(m-1)(4m+1)}. $$
This completes the proof.
\end{proof}

Since $T^2_{k_1,\dots,k_m}$ and
$P^2$ vanish to the same order at $\omega_d$, we have
\begin{multline}\label{tp2}
\lim_{t\rightarrow\omega_d}\frac{T^2_{k_1,\dots,k_m}(t)}{P^2(t)}
=\prod_{i=1}^m  \frac{q^{\frac14(k_i-1)}}{1-q^{\frac 12k_i}}
\prod_{\substack{1\leq i\leq m\\k_i\equiv d\,(2d)}} k_i\\
\times\prod_{\substack{1\leq
  i<j\leq m\\k_i\equiv k_j\,(2d)}}
\left(\frac {k_j-k_i}2\right)^2
\prod_{\substack{1\leq
  i<j\leq m\\k_i\equiv -k_j\,(2d)}}
\left(\frac{k_i+k_j}{2}\right)^2
\prod_{\substack{1\leq i<j\leq 2m\\i\equiv
    j\,(d)}}\frac 1{j-i},
\end{multline}
where, by  the discussion leading to \eqref{pzm},
$$\prod_{\substack{1\leq i<j\leq 2m\\i\equiv
    j\,(d)}}(j-i)=d^{(2m-d)m/d}\prod_{l=1}^{(2m-d)/d}l^{2m-ld}
=d^{(2m-d)m/d}\prod_{l=1}^{(2m-d)/d}l!^d.$$

Plugging \eqref{tp1} and \eqref{tp2} into \eqref{epxo2} gives
\begin{multline*}\sideset{}{'}
\sum_{k_1,\dots,k_m}\,
(-1)^{|\{i;\, d+1\leq k_i\leq 2d-1\, (2d)\}|
}
\prod_{i=1}^m  \frac{q^{\frac14(k_i-1)}}{1-q^{\frac 12k_i}}\\
\times\prod_{\substack{1\leq i\leq m\\k_i\equiv d\,(2d)}} k_i
\prod_{\substack{1\leq
  i<j\leq m\\k_i\equiv k_j\,(2d)}}
\left(\frac {k_j-k_i}2\right)^2
\prod_{\substack{1\leq
  i<j\leq m\\k_i\equiv -k_j\,(2d)}}
\left(\frac{k_i+k_j}{2}\right)^2\\
= (-1)^{\binom{d}{2}\binom{2m/d}{2}}q^{\frac
    14m(m-1)}d^{(2m-d)m/d}m!\prod_{l=1}^{(2m-d)/d}l!^d\,
\triangle(q^{d/2})^{4m^2/d}.
\end{multline*}
Using \eqref{gs} to expand the left-hand side, we arrive at \eqref{mt1}.

Finally we note that, by  Lemma \ref{l1}.b,
\begin{equation}\label{2a}\prod_{\substack{1\leq
  i<j\leq m\\k_i\equiv k_j\,(2d)}}4
\prod_{\substack{1\leq
  i<j\leq m\\k_i\equiv -k_j\,(2d)}}4=\begin{cases}4^{2\binom {m/d}{2}+\frac 12(d-1)\binom{2m/d}{2}}
=2^{m(2m-d-1)/d},& d \text{ odd},\\
4^{\frac 12d\binom{2m/d}{2}}
=2^{m(2m-d)/d},& d \text{ even},
\end{cases}
\end{equation}
which should be used when deriving
\eqref{mt1ab}  from \eqref{mt1}.

{\bf Proof of \eqref{mt2}.}
We  repeat the analysis above, starting with \eqref{opx} 
rather than \eqref{epxo}. 
Consider first the right-hand side of \eqref{opx}. Dividing by
$$  P(t)=\prod_{1\leq i<j\leq 2m+1}(1-t^{j-i}), $$
we have
\begin{equation}\label{vdp}\frac{1}{P(t)}\frac{(q)_\infty^{2m}}{(\sqrt
  q)_\infty^{2m}}\prod_{1\leq i<j\leq 
    2m+1}\frac{\theta(t^{j-i})}{\theta(\sqrt q t^{j-i})}
=\frac{(\sqrt
  q)_\infty}{(q)_\infty}\prod_{i,j=1}^{2m+1}\frac{(qt^{j-i})_\infty}{(\sqrt
  qt^{j-i})_\infty}.\end{equation}

\begin{lemma}
If $d\mid 2m$ or $d\mid 2m+2$, the expression \eqref{vdp} equals
$\triangle(q^{d/2})^{4m(m+1)/d}$ for $t=\omega_d$.
\end{lemma}

\begin{proof}
If $d\mid 2m$, \eqref{vdp} can be written
$$\prod_{i,j=1}^{2m}\frac{(qt^{j-i})_\infty}{(\sqrt
  qt^{j-i})_\infty}\prod_{k=1}^{2m}\frac{(qt^{2m+1-k},qt^{k-2m-1})_\infty}
{(\sqrt qt^{2m+1-k},\sqrt qt^{k-2m-1})_\infty}.
 $$
If $t=\omega_d$, the double product  is computed in  \eqref{rhl} as
$\triangle(q^{d/2})^{4m^2/d}$.
By \eqref{rup}, the single product can be written
$$\prod_{k=1}^{d}\frac{(q\omega_d^{2m+1-k},q\omega_d^{k-2m-1})_\infty^{2m/d}}
{(\sqrt q\omega_d^{2m+1-k},\sqrt q\omega_d^{k-2m-1})_\infty^{2m/d}}
=\frac{(q^d,q^d;q^d)_\infty^{2m/d}}
{(q^{d/2},q^{d/2};q^d)_\infty^{2m/d}}=\triangle(q^{d/2})^{4m/d},$$
which proves the result in this case.

If $d\mid 2m+2$ we may write \eqref{vdp} as
$$\prod_{i,j=1}^{2m+2}\frac{(qt^{j-i})_\infty}{(\sqrt
  qt^{j-i})_\infty}\prod_{k=1}^{2m+2}\frac{(\sqrt qt^{2m+2-k},\sqrt qt^{k-2m-2})_\infty}{(qt^{2m+2-k},qt^{k-2m-2})_\infty}. $$
As above, the double product equals $\triangle(q^{d/2})^{4(m+1)^2/d}$
and the single product  $\triangle(q^{d/2})^{-4(m+1)/d}$, which
  completes the proof.
\end{proof}

\begin{lemma}\label{lm}
 If $d\mid 2m$ or $d\mid 2m+2$, 
 the multiplicity of $t=\omega_d$ as a zero of $P$ is
  $(2m+2-d)m/d$. 
\end{lemma}

\begin{proof}
We must count the pairs $(i,j)$ such that $1\leq i<j\leq 2m+1$ and
$i\equiv j\pod d$. If $j=i+ld$, there are $2m+1-ld$ such pairs, so the 
total number is
 $$\sum_{l=1}^{[2m/d]}(2m+1-ld). $$
 If  $d\mid 2m$,
  we have an arithmetic sum with $2m/d$ terms 
and average value $(2m+2-d)/2$.
If $d\mid 2m+2$, $[2m/d]=(2m+2-d)/d$, unless $d=1$ which is
included in the previous case. This gives
an arithmetic sum with $(2m+2-d)/d$ terms and average value $m$. In
both cases, the result is $(2m+2-d)m/d$.
\end{proof}

We note in passing that the same argument gives,
if $d\mid 2m$,
\begin{subequations}\label{ijp}
\begin{equation}\label{ijp1}\begin{split}\prod_{\substack{1\leq i<j\leq m+1\\i\equiv j\,(d)}}(j-i)
&=d^{(2m+2-d)m/d}\prod_{l=1}^{2m/d}l^{2m+1-ld}\\
&=d^{(2m+2-d)m/d} (2m/d)!
\prod_{l=1}^{(2m-d)/d}l!^d,
\end{split} \end{equation}
while if $d\mid 2m+2$, 
\begin{equation}\label{ijp2}\begin{split}\prod_{\substack{1\leq i<j\leq m+1\\i\equiv j\,(d)}}(j-i)
&=d^{(2m+2-d)m/d}\prod_{l=1}^{(2m+2-d)/d}l^{2m+1-ld}\\
&=d^{(2m+2-d)m/d}
\frac 1{((2m+2-d)/d)!}
\prod_{l=1}^{(2m+2-d)/d}l!^d.
\end{split} \end{equation}
\end{subequations}  
Initially, we only obtain \eqref{ijp2} when $d\neq 1$, but 
when $d=1$ it agrees with \eqref{ijp1} and  thus remains valid.

\begin{lemma}\label{rm}
If $d\mid 2m$ or $d\mid 2m+2$,  the multiplicity of $t=\omega_d$
as a zero of
$$\prod_{i=1}^m(1-t^{k_i})^2(1-t^{2k_i})
\prod_{1\leq i<j\leq
m}(t^{k_j}-t^{k_i})^2(1-t^{k_i+k_j})^2$$
is at least  $(2m+2-d)m/d$, with equality if and only if
\begin{equation}\label{kco}(k_1,\dots,k_m)\simeq (1,2,\dots,m) \mod d.
\end{equation}
\end{lemma}

\begin{proof}
Let $a_i$ be the number of $j$, $1\leq j\leq m$, such that $k_j\equiv
\pm i\pod d$. Here, $i=0,1,2,\dots, (d-1)/2$ or $i=0,1,2,\dots,d/2$
according to whether $d$ is odd or even. If $d$ is odd, the 
multiplicity is 
$$\mu_1=3a_0+4\binom{a_0}{2}+2\binom{a_1}{2}+\dots+2 \binom{a_{(d-1)/2}}{2},
$$
while if $d$ is even, it is
$$\mu_2=3a_0+4\binom{a_0}{2}+2\binom{a_1}{2}+\dots+2 \binom{a_{(d-2)/2}}{2}
+a_{d/2}+4\binom{a_{d/2}}{2}.
$$
In contrast to the case of  Lemma \ref{l1}, the minimum of these
expressions subject to $\sum_ia_i=m$ is achieved at non-integral
values of $a_i$. To prove that the minimum over the natural numbers
equals  $\mu_0=(2m+2-d)m/d$, we 
consider a number of different cases separately. 

If $d$ is odd and $d\mid 2m$ (so that $d\mid m$), we write
$$\mu_1-\mu_0=2\left(a_0-\frac md\right)\left(a_0-\frac md+1\right)+
\sum_{k=1}^{(d-1)/2}\left(a_k-\frac{2m}{d}\right)^2.
 $$
Since $x(x+1)\geq 0$ for integer $x$, we have that $\mu_1\geq \mu_0$
with equality precisely if
$$a_0=m/d,\qquad a_1=\dots=a_{(d-1)/2}=2m/d; $$ 
$a_0=(m-d)/d$ would contradict $\sum a_i=m$.

 If $d$ is odd and $d\mid 2m+2$, we similarly write
$$\mu_1-\mu_0=2\left(a_0-\frac
  {m+1}d\right)\left(a_0-\frac {m+1}d+1\right)+
\sum_{k=1}^{(d-1)/2}\left(a_k-\frac{2(m+1)}{d}\right)^2,$$
which shows that  $\mu_1\geq \mu_0$
with equality precisely if
$$a_0=\frac{m+1}{d}-1,\qquad a_1=\dots=a_{(d-1)/2}=\frac{2(m+1)}d. $$ 

If $d$ is even and  $d\mid 2m$, we write
$$\mu_2-\mu_0=2\left(a_0-\frac md\right)\left(a_0-\frac md+1\right)+
\sum_{k=1}^{(d-2)/2}\left(a_k-\frac{2m}{d}\right)^2+2\left(a_{d/2}-\frac
md\right)^2.
 $$
If $d\mid m$, it follows as before that  $\mu_2\geq \mu_0$
with equality precisely if
$$a_0=a_{d/2}=m/d,\qquad a_1=\dots=a_{(d-2)/2}=2m/d. $$ 
However, if $2m/d$ is odd, we observe that
$$\left(a_0-\frac md\right)\left(a_0-\frac md+1\right)\geq -\frac 14,\qquad
\left(a_{d/2}-\frac
md\right)^2\geq \frac 14, $$
so we still have $\mu_2\geq \mu_0$, but with equality for
$$a_0=\frac md-\frac12,\qquad a_1=\dots=a_{(d-2)/2}=\frac{2m}d,\qquad
 a_{d/2}=\frac md+\frac12. $$ 

Finally, if $d$ is an even divisor of $2m+2$,
\begin{equation*}\begin{split}\mu_2-\mu_0&=2\left(a_0-\frac {m+1}d\right)\left(a_0-\frac {m+1}d+1\right)+
\sum_{k=1}^{(d-2)/2}\left(a_k-\frac{2(m+1)}{d}\right)^2\\
&\quad+2\left(a_{d/2}-\frac
{m+1}d\right)^2,
 \end{split}\end{equation*}
and we conclude as above that $\mu_2\geq \mu_0$ with equality precisely when
$$a_0=\frac {m+1}d-1,\qquad a_1=\dots=a_{(d-2)/2}=\frac{2(m+1)}d,\qquad
a_{d/2}=\frac{m+1}d $$ 
if $(2m+2)/d$ is even and when
$$a_0=a_{d/2}=\frac {m+1}d-\frac 12,\qquad a_1=\dots=a_{(d-2)/2}=\frac{2(m+1)}d $$ 
if  $(2m+2)/d$ is odd. 

In each case, the desired result now  follows from Lemma \ref{l2}.b.
\end{proof}

We conclude from Lemma \ref{lm} and Lemma \ref{rm} that
\begin{multline}\label{os}
\frac{1}{m!}\,\sideset{}{'}
\sum_{k_1,\dots,k_m}\,\frac{T^1_{k_1,\dots,k_m}(\omega_d)}{P^1(\omega_d)}
\lim_{t\rightarrow\omega_d}\frac{T^2_{k_1,\dots,k_m}(t)}{P^2(t)}\\
=  q^{\frac
    14m(m+1)}\omega_d^{-\frac13m(m+1)(2m+1)}\triangle(q^{d/2})^{4m(m+1)/d}
,\end{multline}
where the sum is over positive integers satisfying \eqref{kco},
and where
{\allowdisplaybreaks
\begin{align*}
T^1&=\prod_{i=1}^m t^{-2mk_i}
\prod_{\substack{1\leq i<j\leq
    m\\k_i\equiv k_j\,(d)}}t^{2k_i}
\prod_{\substack{1\leq i\leq m\\k_i\not\equiv 0\,(d)}}
(1-t^{k_i})^2\prod_{\substack{1\leq i\leq m\\2k_i\not\equiv 0\,(d)}}
(1-t^{2k_i})
\\*
&\quad\times\prod_{\substack{1\leq
  i<j\leq m\\k_i\not\equiv k_j\,(d)}}
(t^{k_i}-t^{k_j})^2
\prod_{\substack{1\leq
  i<j\leq m\\k_i\not\equiv -k_j\,(d)}}
(1-t^{k_i+k_j})^2,\\
T^2&= \prod_{i=1}^m  \frac{q^{\frac12 k_i}}{1-q^{k_i}}
\prod_{\substack{1\leq i\leq m\\k_i\equiv 0\,(d)}}
(1-t^{k_i})^2
\prod_{\substack{1\leq i\leq m\\2k_i\equiv 0\,(d)}}
(1-t^{2k_i})
\\*
&\quad\times\prod_{\substack{1\leq
  i<j\leq m\\k_i\equiv k_j\,(d)}}
(1-t^{k_j-k_i})^2
\prod_{\substack{1\leq
  i<j\leq m\\k_i\equiv -k_j\,(d)}}
(1-t^{k_i+k_j})^2,\\
P^1(t)&=\prod_{\substack{1\leq i<j\leq 2m+1\\i\not\equiv
    j\,(d)}}(1-t^{j-i}),\\
P^2(t)&=\prod_{\substack{1\leq i<j\leq 2m+1\\i\equiv
    j\,(d)}}(1-t^{j-i}).
\end{align*}
}

The expression $T^1_{k_1,\dots,k_m}(\omega_d)$ is visibly symmetric in
the $k_i$ and invariant under $k_i\mapsto k_i+d$. It is
straight-forward to check that it is odd in $k_1$ if $2k_1\not\equiv 0\pod
d$. Thus, using also Lemma \ref{l2}.c, we have
$$T^1_{k_1,\dots,k_m}(\omega_d)=(-1)^{|\{i;\,(d+1)/2\leq k_i\leq
  d-1\,(d)\}|
+\binom{d-1}2\binom{2M/d}{2}}
T^1_{1,2,\dots,m}(\omega_d), $$
where $M=m$ if $d\mid 2m$ and $M=m+1$ if  $d\mid 2m+2$.
Moreover, similarly as in Lemma \ref{tpl}, it is straight-forward to check
  that
$$\frac{T^1_{1,2,\dots,m}(\omega_d)}{P^1(t)}=\prod_{i=1}^mt^{-2mi}\prod_{1\leq
i<j\leq m}t^{2i}=t^{-\frac 13m(m+1)(2m+1)}. $$
Thus, 
$$
\frac{T^1_{k_1,\dots,k_m}(\omega_d)}{P^1(\omega_d)}
=(-1)^{|\{i;\,(d+1)/2\leq k_i\leq
  d-1\,(d)\}|
+\binom{d-1}2\binom{2M/d}{2}}\omega_d^{-\frac 13m(m+1)(2m+1)}.
$$

We also have
\begin{equation*}\begin{split}
\lim_{t\rightarrow\omega_d}\frac{T^2_{k_1,\dots,k_m}(t)}{P^2(t)}
&= \prod_{i=1}^m  \frac{q^{\frac12 k_i}}{1-q^{k_i}}
\prod_{\substack{1\leq i\leq m\\k_i\equiv 0\,(d)}}2k_i^3
\prod_{\substack{1\leq i\leq m\\k_i\equiv d/2\,(d)}}2k_i
\\
&\quad\times\prod_{\substack{1\leq
  i<j\leq m\\k_i\equiv k_j\,(d)}}
(k_j-k_i)^2
\prod_{\substack{1\leq
  i<j\leq m\\k_i\equiv -k_j\,(d)}}
(k_i+k_j)^2\prod_{\substack{1\leq i<j\leq 2m+1\\i\equiv
    j\,(d)}}\frac 1{j-i},
\end{split}\end{equation*}
where the final factor is given by \eqref{ijp}.

Plugging all this into \eqref{os}, we obtain
\begin{multline*}
\sideset{}{'}
\sum_{k_1,\dots,k_m}\,
(-1)^{\left|\left\{i;\,(d+1)/2\leq k_i\leq
  d-1\,(d)\right\}\right|}
\prod_{i=1}^m  \frac{q^{\frac12 k_i}}{1-q^{k_i}}
\prod_{\substack{1\leq i\leq m\\k_i\equiv 0\,(d)}}2k_i^3
\prod_{\substack{1\leq i\leq m\\k_i\equiv d/2\,(d)}}2k_i
\\
\times\prod_{\substack{1\leq
  i<j\leq m\\k_i\equiv k_j\,(d)}}
(k_j-k_i)^2
\prod_{\substack{1\leq
  i<j\leq m\\k_i\equiv -k_j\,(d)}}
(k_i+k_j)^2,
\\
=Cd^{(2m+2-d)m/d} q^{\frac
    14m(m+1)}
m!\,
\triangle(q^{d/2})^{4m(m+1)/d}
,\end{multline*}
with $C$ as in \eqref{c}. Expanding the left-hand side using
\eqref{ope}, we  arrive at \eqref{mt2}.

Finally, we note that, by Lemma \ref{l2}.b,
\begin{equation}\label{2b}
\prod_{\substack{1\leq i\leq m\\k_i\equiv 0\,(d)}}2
\prod_{\substack{1\leq i\leq m\\k_i\equiv d/2\,(d)}}2
=\begin{cases}2^{2m/d}, & d \text{ even},\ d\mid 2m\\
2^{m/d}, & d \text{ odd},\ d\mid 2m\\
2^{(2m+2-d)/d}, & d \text{ even},\ d\mid 2m+2\\
2^{(m+1-d)/d}, & d \text{ odd},\ d\mid 2m+2,
\end{cases}
\end{equation}
which should be used in deriving \eqref{mt2ab} from
\eqref{mt2}.

\end{document}